\documentclass[11pt]{article}
\input epsf.tex
\usepackage{amssymb,latexsym,times}
\catcode`\@=11
\renewcommand\subsection{\@startsection{subsection}{2}{\z@}%
                                     {-3.25ex\@plus -1ex \@minus -.2ex}%
                                     {-0.01 mm}
                                     {\normalfont\large\bfseries}}

\catcode`\@=11
\renewcommand\subsubsection{\@startsection{subsubsection}{2}{\z@}%
                                     {-3.25ex\@plus -1ex \@minus -.2ex}%
                                     {-0.01 mm}
                                     {\normalfont\bfseries}}

\newtheorem{example}{Example}
\newtheorem{definition}[example]{Definition}
\newtheorem{remark}{Remark}
\newtheorem{conjecture}{Conjecture}

\def\adots{\mathinner{\mkern2mu\raise1pt\hbox{.}
\mkern3mu\raise4pt\hbox{.}\mkern1mu\raise7pt\hbox{.}}}
\def\<{\langle\,}
\def\>{\,\rangle}

\def\shuf{\sqcup\!\sqcup}
\def\ie{{\em i.e. }}

\def\SG{\mathfrak S}

\def\a{\alpha}

\def\b{\beta}

\def\N{{\mathbb N}}
\def\Z{{\mathbb Z}}
\def\C{{\mathbb C}}

\def\Q{{\mathbb Q}}

\def\F{{\cal F}}
\def\P{{\cal P}}

\def\B{{\bf B}}
\def\Bs{{\bf B}^*}

\def\LL{{\mathfrak L}}

\def\g{\mathfrak g}

\def\Sl{\mathfrak{sl}}

\def\n{\mathfrak n}

\def\slchap{\widehat{\mathfrak{sl}}}

\def\H{\widehat H}

\def\LL{{\cal L}}
\def\L{\Lambda}

\def\<{\langle}
\def\>{\rangle}

\def\BB{{\cal B}}
\def\CC{{\cal C}}
\def\m{{\bf m}}

\def\p{{\bf p}}

\def\b{{\bf b}}

\def\RR{{\cal R}}

\def\le{\leqslant}
\def\ge{\geqslant}

\def\c{{\mathfrak c}}
\def\shuf{*}

\title{\bf Imaginary vectors in the dual canonical basis of $U_q(\n)$} 

\author{Bernard {\sc Leclerc}}

\date{}

\begin{document}
\maketitle


\begin{abstract}\noindent
Let $\n$ be the maximal nilpotent subalgebra of a simple
complex Lie algebra $\g$. 
We introduce the notion of imaginary vector in the dual canonical
basis of $U_q(\n)$, and we give examples of such vectors for
types $A_n (n\ge 5)$, $B_n (n\ge 3)$, $C_n (n\ge 3)$,
$D_n (n\ge 4)$, and all exceptional types.
This disproves a conjecture of Berenstein and Zelevinsky
about $q$-commuting products of vectors of the dual canonical basis.
It also shows the existence of finite-dimensional irreducible
representations of quantum affine algebras whose tensor square
is not irreducible.
\end{abstract}

\vskip 0.6cm

\section{Introduction} 
\label{SECT1}

Let $\g$ be a simple complex Lie algebra, $U_q(\g)$ 
the corresponding quantized enveloping algebra,
$\n$~a maximal nilpotent subalgebra of $\g$, and $U_q(\n)$
the corresponding subalgebra of $U_q(\g)$.  
Let ${\bf B}$ be the canonical basis of $U_q(\n)$ \cite{L,K0},  
and let $\Bs$ be the basis dual to ${\bf B}$ for the natural scalar
product on $U_q(\n)$.
In this note, we discuss the multiplicative structure of~$\Bs$.
We write $q^\Z\Bs = \{q^mb\ |\ m\in\Z,\ b\in \Bs\}$.
For $u, v$ in $U_q(\n)$ we also use the short-hand notation 
$u \cong v$ to mean that there exists $m\in\Z$ such that $u = q^mv$.

\begin{definition} {\rm 
We say that $b\in\Bs$ is {\em real} if $b^2\in q^\Z\Bs$. 
Otherwise, we say that $b$ is {\em imaginary}. 
} 
\end{definition}
In \cite{BZ} it was conjectured that, given $b_1, b_2 \in \Bs$,
the product $b_1b_2$ is in $q^\Z\Bs$ if and only if
$b_1$ and $b_2$ $q$-commute with each other, \ie
$b_2b_1 \cong b_1b_2$.
In particular, taking $b_1=b_2$, this conjecture would imply that 
all vectors of $\Bs$ are real.
In Section~\ref{SECT2} below, 
we give examples of imaginary vectors in $U_q(\n)$ for $\g$
of type $A_5$, $B_3$, $C_3$, $D_4$ and $G_2$.
Using the natural embeddings $U_q(A_5) \subset U_q(A_n) \ (n\ge 6)$,
$U_q(B_3) \subset U_q(B_n) \ (n\ge 4)$,
$U_q(C_3) \subset U_q(C_n) \ (n\ge 4)$,
$U_q(D_4) \subset U_q(D_n) \ (n\ge 5)$,
$U_q(D_4) \subset U_q(E_n) \ (n=6,7,8)$,
$U_q(B_3) \subset U_q(F_4)$,
it follows that there exist imaginary vectors for all types
except $A_1,A_2,A_3,A_4,B_2$.

For types $A_3,B_2$ the basis $\Bs$ is explicitely described
in \cite{BZ} and \cite{RZ,Z,C} respectively, and it consists in all 
$q$-commuting products of elements of a certain finite subset $\P$ of 
$\Bs$.
The elements of $\P$ may be regarded as the ``prime vectors'' of $\Bs$.
A similar description also exists in type $A_4$ \cite{Z2}. 
From these descriptions one sees that all vectors of $\Bs$
are real if $\g$ is of type $A_n (n\le 4)$ or $B_2$.

In Section~\ref{SECT3}, we propose a conjecture regarding
the multiplication by a real element of~$\Bs$.
This conjecture would imply that the Berenstein-Zelevinsky
conjecture is true when one of the two factors is real.
It would also yield an analogue of the Kashiwara crystal
graph operator $\widetilde{e}_i$ for any real vector of $\Bs$.

Our original motivation for investigating the multiplicative
properties of $\Bs$ came from some questions in
the theory of finite dimensional representations of
quantum affine algebras.
Let $\widehat{\g}$ be an affine Lie algebra and $U_q(\widehat{\g})$
its quantum enveloping algebra. 
Let $V, W$ be finite dimensional irreducible representations of $U_q(\widehat{\g})$.
It is known that if $V\otimes W$ is irreducible then it is isomorphic
to $W\otimes V$.
Using our example of imaginary vector in type $A_n$, we 
are able to show that the converse is not true for 
$\widehat{\g}=\widehat{\Sl}_N$, namely, we can exhibit an 
explicit irreducible representation $V$ such that
$V\otimes V$ is not irreducible (see Section~\ref{SECT4}).

The name ``imaginary vector'' has been
inspired by the recent seminal work of Fomin and Zelevinsky on 
cluster algebras \cite{FZ,Z3}.
According to Zelevinsky \cite{Z2}, the general picture
should be as follows.
Let $G$ be the simply connected complex simple Lie group
with Lie algebra $\g$ and $N$ its maximal unipotent subgroup
with Lie algebra $\n$.
The specialization $\C[N]$ of $U_q(\n)$ at $q=1$ (using the
$\Z[q,q^{-1}]$-lattice $\LL^*$ spanned by $\B^*$) is expected to be
a cluster algebra, and to have finite cluster type only
for $\g$ of type $A_n (n\le 4)$ or $B_2$.
In this ``finite case'' the prime elements of $\Bs$ not belonging
to the $q$-center are naturally labelled
by the positive roots and the negatives of the simple roots
of a certain complex simple Lie algebra~$\c$.
In the ``infinite case'', it is expected that the imaginary prime vectors
will correspond to imaginary roots of some infinite-dimensional
Kac-Moody algebras.

\section{Examples of imaginary vectors} 
\label{SECT2}   

\subsection{}
Let $e_i \ (i=1,\ldots ,r)$ be the Chevalley generators of $U_q(\n)$.
Let $w_0$ denote the longest element of the Weyl group $W$
of $\g$. To a reduced decomposition $w_0=s_{i_1}\cdots s_{i_n}$
is associated a convex ordering $\beta_1<\cdots <\beta_n$
of the positive roots of $\g$ and a PBW-type basis of $U_q(\n)$ \cite{L}:
\[
E(\m)=E(\beta_1)^{(m_1)}\cdots E(\beta_n)^{(m_n)},
\qquad (\m=(m_1,\ldots ,m_n)\in \N^n).
\]
Let $\{E^*(\m)\}$ be the dual PBW-basis, and let $\b(\m)$ denote
the unique element of $\Bs$ such that $\b(\m)-E^*(\m)\in qL^*$,
where $L^*$ is the $\Z[q]$-lattice spanned by $\Bs$.

We shall use repeatedly the following known fact (see \cite{Rei}). 
For any $\m, \p \in \N^n$, the $\Bs$- expansion of $\b(\m)\b(\p)$
contains the vector $\b(\m+\p)$ with coefficient a power of $q$.
Hence, if $\b(\m)\b(\p)\in q^{\Z}\Bs$ then
$\b(\m)\b(\p)\cong\b(\m+\p)$.

For $k\in\N^*$ and $b=\b(\m)\in\Bs$ we set
$k\m = (km_1,\ldots ,km_n)$ and $b^{[k]}=\b(k\m)$.
Thus $b$ is real if and only if $b^2 \cong b^{[2]}$. 

\subsection{}
The following examples of imaginary vectors have been calculated with {\tt maple},
using the algorithms described in \cite{Le} for calculating
the elements of $\Bs$ and their products.
In \cite{Le}, the basis $\Bs$ is studied in the realization
of $U_q(\n)$ in terms of quantum shuffles, following 
Rosso \cite{R1,R2} and Green \cite{G}.
We will now explain very briefly the main features of
this approach, refering the interested reader to \cite{Le}
for a detailed exposition.

In the quantum shuffle realization a vector $v\in U_q(\n)$ is represented 
by a $\Q(q)$-linear combination of words $w_{i_1}\ldots w_{i_k}$ over
the alphabet $\{w_1,\ldots ,w_r\}$.
In particular, $e_i$ is identified with the letter $w_i$.
In the sequel we shall write $w[i_1,\ldots ,i_k]$
rather than $w_{i_1}\cdots w_{i_k}$.
The multiplication of $U_q(\n)$ translates into the
bilinear operation $\shuf$ defined on words by
\begin{equation}\label{eq1}
    w[i_1,\ldots ,i_m]\shuf w[i_{m+1},\ldots ,i_{m+n}]
=
\sum_\sigma q^{-e(\sigma)} w[i_{\sigma(1)},\ldots ,i_{\sigma(m+n)}] 
\end{equation}
where 
the sum runs over the $\sigma \in \SG_{m+n}$ such that
$\sigma(1)<\cdots < \sigma(m)$ and $\sigma(n+1)<\cdots < \sigma(m+n)$,
and 
\[
 e(\sigma) = \sum_{k\le m<l;\ \sigma(k)<\sigma(l)}
(\a_{i_{\sigma(k)}}\,,\a_{i_{\sigma(l)}}).
\]
Here, $\a_1,\ldots ,\a_r$ are the simple roots of $\g$ and
$(\cdot , \cdot)$ denotes the bilinear form on the root lattice
such that the entries of the Cartan matrix of $\g$ are
given by
\[
a_{ij} = {2(\a_i\,,\a_j)\over (\a_i\,,\a_i)} = {(\a_i\,,\a_j)\over
  d_i}, \qquad (1\le i,j \le r)
\]
where $d_i = (\a_i\,,\a_i)/2 \in \{1,2,3\}$.
When $q=1$, Equation (\ref{eq1}) is nothing else than the classical
shuffle product \cite{Lo,Reu}.

Following Rosso \cite{R3} and Lalonde and Ram \cite{LR}, 
it is shown in \cite{Le} that to each total ordering of the
set of simple roots corresponds a certain (dual) PBW basis,
whose root vectors can be readily calculated using (\ref{eq1}) by taking
some iterated $q$-bracketings.
(In the calculations below we take the natural order
$\a_1<\cdots<\a_r$, which determines our choice of PBW basis.)
Then, the characterization of $\Bs$ given in \cite{Le} (Proposition~41)
yields an easy algorithm for calculating it.
Moreover, in this setting, the products of elements
of $\Bs$ and the calculation of their $\Bs$-expansion 
are rather simple to compute, thanks to the elementary formula
(\ref{eq1}) and to the fact that the words $w[i_1,\ldots ,i_k]$
are linearly independent.

We will indicate below the main steps of the calculation of 
an imaginary vector of type $G_2$ in the quantum shuffle realization.
In that case, the expressions of the vectors are so small that
the computation could be done by hand.
For the other types, the calculations are very similar but
much bigger, so we shall omit them and give only the result.

\subsection{}\label{TG2}
{\it type $G_2$:}
The numbering of the simple roots $\alpha_1, \alpha_2$ is shown by the
following Dynkin diagram:
\begin{center}
\leavevmode
\epsfxsize =2.3cm
\epsffile{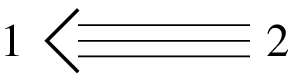}
\end{center}
Thus, we have $(\alpha_1,\alpha_1)=2,\ (\alpha_2,\alpha_2)=6$. 
We choose the reduced decomposition 
\[
w_0=s_1s_2s_1s_2s_1s_2,
\]
which gives the following ordering of the positive roots:
\[
\alpha_1<3\alpha_1+\alpha_2<2\alpha_1+\alpha_2<3\alpha_1+2\alpha_2
<\alpha_1+\alpha_2<\alpha_2.
\]
Let $b=\b(1,0,0,0,1,0)$.
In the quantum shuffle realization, we have 
\[
b=w[1,2,1],
\]
hence, using Equation (\ref{eq1}) we obtain
\[
b^2=
2(1+q^{-2})w[1,2,1,1,2,1]+(q^4+2q^2+1+q^{-2}+2q^{-4}+q^{-6})w[1,1,2,2,1,1]
\]
\[
+\ (q+2q^{-1}+q^{-3})w[1,2,1,2,1,1]+(q+2q^{-1}+q^{-3})w[1,1,2,1,2,1].
\]
On the other hand, using the algorithm mentioned above for calculating
$\Bs$, we get
\[
b^{[2]}=
(q+q^{-1})w[1,2,1,1,2,1]+(q^5+2q^3+q+q^{-1}+2q^{-3}+q^{-5})w[1,1,2,2,1,1]
\]
\[
+\ (q^2+2+q^{-2})w[1,2,1,2,1,1]+(q^2+2+q^{-2})w[1,1,2,1,2,1].
\]
Therefore, $b^2 \not\cong b^{[2]}$. More precisely, putting  
\[
z = (q+q^{-1})w[1,2,1,1,2,1] =\b(1,0,1,0,1,0) \,,
\]
we obtain
\[
b^2 = q^{-1}(b^{[2]}+z),
\]
hence $b$ is imaginary.

\subsection{}
{\it type $B_3$:}
The numbering of the simple roots is given by the following Dynkin diagram:
\begin{center}
\leavevmode
\epsfxsize =3.8cm
\epsffile{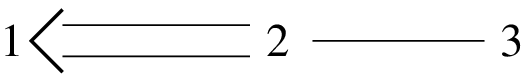}
\end{center}
Thus, we have 
\[
(\alpha_1,\alpha_1)=2,\ (\alpha_2,\alpha_2)=(\alpha_3,\alpha_3)=4.
\] 
We choose the reduced decomposition $w_0=s_1s_2s_3s_1s_2s_1s_3s_2s_3$,
which gives the following ordering of the positive roots:
\[
\alpha_1<2\alpha_1+\alpha_2<2\alpha_1+\alpha_2+\alpha_3<\alpha_1+\alpha_2<
2\alpha_1+2\alpha_2+\alpha_3<\alpha_1+\alpha_2+\alpha_3<\alpha_2
<\alpha_2+\alpha_3<\alpha_3.
\]
Let $b=\b(0,1,0,0,0,0,0,1,0)$, $z=\b(0,1,0,0,1,0,0,1,0)$.
We have 
\[
b^2 = q^{-2}(b^{[2]}+z),
\]
hence $b$ is imaginary.

\subsection{}
{\it type $C_3$:}
The numbering of the simple roots is given by the following Dynkin diagram:
\begin{center}
\leavevmode
\epsfxsize =3.8cm
\epsffile{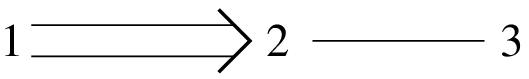}
\end{center}
Thus, we have 
\[
(\alpha_1,\alpha_1)=4,\ (\alpha_2,\alpha_2)=(\alpha_3,\alpha_3)=2.
\] 
We choose the reduced decomposition $w_0=s_1s_2s_1s_3s_2s_1s_3s_2s_3$,
which gives the following ordering of the positive roots:
\[
\alpha_1<\alpha_1+\alpha_2<\alpha_1+2\alpha_2<\alpha_1+\alpha_2+\alpha_3<
\alpha_1+2\alpha_2+\alpha_3<\alpha_1+2\alpha_2+2\alpha_3<\alpha_2
<\alpha_2+\alpha_3<\alpha_3.
\]
Let $b=\b(0,1,0,0,0,0,0,1,0)$, $z=\b(0,1,0,0,1,0,0,1,0)$.
We have 
\[
b^2 = q^{-1}(b^{[2]}+z),
\]
hence $b$ is imaginary.

\subsection{}
{\it type $D_4$:}
The numbering of the simple roots is given by the following Dynkin diagram:
\begin{center}
\leavevmode
\epsfxsize =3.2cm
\epsffile{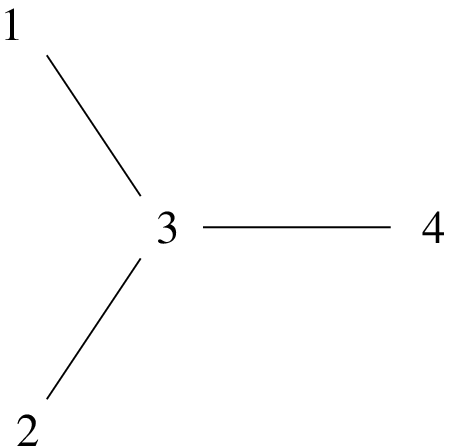}
\end{center}
We choose the reduced decomposition 
$
w_0=s_1s_3s_2s_4s_3s_1s_4s_3s_2s_4s_3s_4,
$
which gives the following ordering of the positive roots:
\[
\alpha_1<\alpha_1+\alpha_3<\alpha_1+\alpha_2+\alpha_3<\alpha_1+\alpha_3+\alpha_4<
\alpha_1+\alpha_2+\alpha_3+\alpha_4\qquad\qquad
\]
\[
\qquad\qquad
<\alpha_1+\alpha_2+2\alpha_3+\alpha_4<
\alpha_2<\alpha_2+\alpha_3<\alpha_2+\alpha_3+\alpha_4<\alpha_3
<\alpha_3+\alpha_4<\alpha_4.
\]
Let 
$b=\b(0,1,0,0,0,0,1,0,0,0,1,0),\
z=\b(0,1,0,0,1,0,0,1,0,0,1,0).
$
We have 
\[
b^2 = q^{-1}(b^{[2]}+z),
\]
hence $b$ is imaginary.

\subsection{}\label{A5}
{\it type $A_5$:}
The simple roots are denoted 
$\alpha_1, \alpha_2, \alpha_3,\alpha_4,\alpha_5$,
and numbered in the natural way.
We choose the reduced decomposition 
\[
w_0=s_1s_2s_3s_4s_5s_1s_2s_3s_4s_1s_2s_3s_1s_2s_1,
\]
which gives the following ordering of the positive roots:
\[
\alpha_1<\alpha_1+\alpha_2<\alpha_1+\alpha_2+\alpha_3<\alpha_1+\alpha_2+\alpha_3+\alpha_4<
\alpha_1+\alpha_2+\alpha_3+\alpha_4+\alpha_5<
\alpha_2<\alpha_2+\alpha_3<\qquad
\]
\[
\qquad
\alpha_2+\alpha_3+\alpha_4<\alpha_2+\alpha_3+\alpha_4+\alpha_5<
\alpha_3<\alpha_3+\alpha_4<\alpha_3+\alpha_4+\alpha_5<
\alpha_4<\alpha_4+\alpha_5<\alpha_5
\]
Let 
\begin{eqnarray*}
b&=&\b(0,1,0,0,0,0,0,1,0,1,0,0,0,1,0),\\
z&=&\b(0,1,0,1,0,0,1,0,1,0,1,0,0,1,0).
\end{eqnarray*}
We have 
\[
b^2 = q^{-2}(b^{[2]}+z),
\]
hence $b$ is imaginary.
\begin{remark}
{\rm
Let $\B$ be the canonical basis for $\g$ of type $A_5$,
and let $c$ and $t$ be the elements of $\B$ such that
$(c,b^{[2]})=(t,z)=1$.
In \cite{KS}, Kashiwara and Saito prove that $t$ is a 
counterexample to a problem of Lusztig \cite{L2}, namely
they show that the singular support $SS(L_t)$ of the irreducible
perverse sheave $L_t$ corresponding to $t$ is not irreducible. 
More precisely, they show that
$$SS(L_t) = \L_t \cup \L_c,$$ where $\L_t$ and $\L_c$
are the irreducible components of the nilpotent variety
associated to $t$ and $c$.
It would be interesting to understand the connection between 
Lusztig's problem and the conjecture of Berenstein and Zelevinsky.
}
\end{remark}
\begin{remark}
{\rm
In all the previous examples the vector $z\in\Bs$
occuring in the expansion of the imaginary vector $b$ is in fact
an element of the $q$-center of $U_q(\n)$, that is, $ze_i \cong e_iz$
for every $i\in \{1,\ldots ,r\}$
(see \cite{C} for a description of the $q$-center).
This needs not to be true in general, since for example a $q$-central element
in $U_q(\n)$ of type $A_5$ will no longer be $q$-central when we
embed $U_q(\n)$ in $U_q(\n')$ of type $A_6$. 
}
\end{remark}

\section{Real vectors and generalized crystal operators} 
\label{SECT3}   

\subsection{} Let $b_1,b_2$ be elements of $\Bs$. Suppose that
$b_1$ is real and $b_1b_2\not\in q^\Z\Bs$.
\begin{conjecture}\label{CONJ1}
The expansion of $b_1b_2$ on $\Bs$ is of the form
\[
b_1b_2=q^{m}\,b'+q^{s}\,b''
       +\sum_{c\not = b',b''} \gamma_{b_1b_2}^{c}(q) \,c,
\]
where $b'\not =b''$, $m,s\in\Z$, $m<s$, $\gamma_{b_1b_2}^{c}(q)\in\Z[q,q^{-1}]$, and
for all $c\in\Bs$ such that $\gamma_{b_1b_2}^{c}(q)\not = 0$
\[ 
\gamma_{b_1b_2}^{c}(q)\in q^{m+1}\Z[q]\cap q^{s-1}\Z[q^{-1}].
\]
\end{conjecture} 
It is known that the $\Bs$-expansion of $b_2b_1$ can be obtained
from that of $b_1b_2$ by changing $q$ into $q^{-1}$ (and rescaling
by an appropriate power of $q$) \cite{Rei}.
If $b_1b_2\not\in q^\Z\Bs$ and Conjecture~\ref{CONJ1} is true, 
we see that $b_2b_1$ cannot be proportional to $b_1b_2$ since
the term with the highest power of $q$ in $b_1b_2$ is $b''$,
while that in $b_2b_1$ is $b'\not = b''$. 
Hence Conjecture~\ref{CONJ1} would imply that the
Berenstein-Zelevinsky conjecture holds true when one of the two vectors
is real.
It is proved in \cite{LNT} that if $\g$ is of type $A_n$, the 
quantum flag minors are real elements of $\Bs$.
Hence our conjecture agrees with the recent result of Caldero
\cite{C2} showing that if a quantum flag minor $q$-commutes with an arbitrary element
of $\Bs$, then their product belongs to $q^\Z\Bs$.

\begin{figure}[t]
\begin{center}
\leavevmode
\epsfxsize =11cm
\epsffile{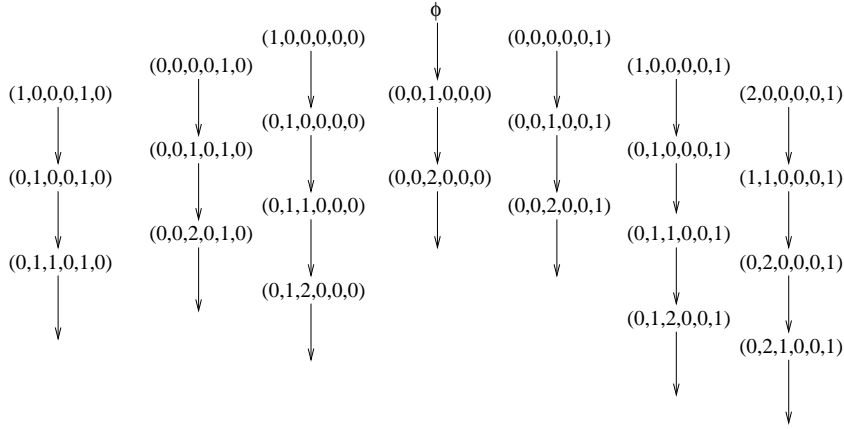}
\end{center}
\caption{\label{FIG1} The first $b_1$-strings with
$b_1=\b(0,0,1,0,0,0)$ in type $G_2$.}
\end{figure}
\subsection{}\label{SS3.2}
Using results of Kashiwara (\cite{K}, 5.3.1), it is easy to check
Conjecture~\ref{CONJ1} in the case where $b_1=e_i$
is a Chevalley generator of $U_q(\n)$ 
(note that $e_i$ is always a real element of $\Bs$).
In that case, 
$b'=\widetilde{e}_i(b_2)$, 
where $\widetilde{e}_i$ is the Kashiwara operator associated with
the simple root vector~$e_i$.
(Here we are abusing notation by denoting in the same way a vector
of $\Bs$ and the corresponding element of the crystal basis at $q=0$.)
Similarly
$b''=\widetilde{e}^{\,\tau}_i(b_2)$, 
where $\widetilde{e}_i^{\,\tau}=\tau\ \widetilde{e}_i\tau$ and $\tau$ is the
involution
of $\Bs$ induced by the anti-automorphism 
of $U_q(\n)$ which fixes the $e_i$'s.

Assuming Conjecture~\ref{CONJ1}, let us write $b'=b_1\diamond b_2$
and $b''=b_2\diamond b_1$.
(If $b_1b_2\cong c\in\Bs$, we put $b_1\diamond b_2=b_2\diamond b_1=c$.) 
For a fixed real $b_1$, we may regard the
maps $b \mapsto b_1\diamond b$ and $b \mapsto b\diamond b_1$ as
generalized left and right Kashiwara operators associated to $b_1$.
(Note that for every positive root $\beta$, the vector $E^*(\beta)$
of the dual PBW-basis is a real vector in $\Bs$.
Thus, in particular, we obtain conjectural analogues of Kashiwara operators 
for all positive roots).

\begin{conjecture}\label{CONJ2}
The maps $b \mapsto b_1\diamond b$ and $b \mapsto b\diamond b_1$ are injective.
\end{conjecture}
Conjecture~\ref{CONJ2} implies that we can partition the set $\Bs$
into infinite $b_1$-strings:
\[
b \longrightarrow b_1\diamond b \longrightarrow b_1 \diamond (b_1\diamond b)
\longrightarrow \cdots
\]
with $b$ running over $\Bs - (b_1 \diamond \Bs)$.
This is illustrated in Figure~\ref{FIG1} for $b_1=\b(0,0,1,0,0,0)$ in
type~$G_2$. We also have a similar string decomposition corresponding to the
right operation $b\mapsto b\diamond b_1$.

\subsection{}
These two conjectures are based on extensive computer calculations.

In particular, Conjecture~\ref{CONJ1} has been checked systematically 
in type $G_2$ for all vectors $b_1$, $b_2$ of principal 
degree~$\le 7$ which are not both imaginary.
Note that among the $116$ vectors of $\Bs$ of degree~$\le 7$,
there are $10$ imaginary vectors:
\begin{eqnarray*}
&&\b(1,0,0,0,1,j),\quad (0\le j\le 4),\\
&&\b(2,0,0,0,2,j),\quad (0\le j\le 1),\\
&&\b(j,1,0,0,0,1),\quad (0\le j\le2).
\end{eqnarray*}
Only $3$ of them are prime, namely $\b(1,0,0,0,1,0)$,
$\b(2,0,0,0,2,0)$ and $\b(0,1,0,0,0,1)$.

\section{Applications to affine Hecke algebras and quantum affine algebras}
\label{SECT4}

One important motivation for studying the multiplicative properties
of $\Bs$ comes from the representation theory of affine Hecke
algebras and quantum affine algebras \cite{LNT}.

\subsection{}
Let $H_m(t)$ be the affine Hecke algebra of type $GL_m$ over $\C$, 
where $t\in\C^*$ is of infinite multiplicative order.
It has invertible generators $T_1,\ldots ,T_{m-1},y_1,\ldots ,y_m$
subject to the relations
\begin{eqnarray*}
&&T_iT_{i+1}T_i=T_{i+1}T_iT_{i+1},\hskip 1.1cm (1\le i\le m-2),\label{EQ_T1}\\
&&T_iT_j=T_jT_i,\hskip 2.9cm (\vert i-j\vert>1),\label{EQ_T2}\\
&&(T_i-t)(T_i+1)=0, \hskip 1.6cm(1\le i\le m-1),\label{EQ_T3}\\
&&y_iy_j=y_jy_i, \hskip 3.1cm (1\le i,j\le m),\label{EQ_YY}\\
&&y_jT_i=T_iy_j, \hskip 3cm (j \not= i,i+1),\\
&&T_iy_iT_i=t\,y_{i+1}, \hskip 2.5cm (1\le i\le m-1).\label{EQ_Y}
\end{eqnarray*}
Let $\CC_m$ denote the category of finite-dimensional
$H_m(t)$-modules, and $\CC_{m,r}$ the sub-category consisting 
of those modules for which all eigenvalues of the pairwise commutative 
generators $y_1,\ldots ,y_m$ belong to $\{t,t^2,\ldots ,t^r\}$.
Let $\RR_r$ be the direct sum of the complexified Grothendieck groups
of the $\CC_{m,r} \ (m\in\N)$, endowed with the multiplication given by the 
induction functors $\CC_{m,r} \times \CC_{n,r} \longrightarrow \CC_{m+n,r}$.
The algebra $\RR_r$ is graded by deciding that the classes of the
$H_m(t)$-modules have degree $m$.
 
Let $N$ be the group of upper triangular unipotent $(r+1)\times (r+1)$-matrices,
$\n$ the Lie algebra of $N$, and $U_q(\n)$ the
corresponding quantum enveloping algebra.
(This is the positive part of a quantum enveloping algebra of type
$A_r$.)
Let $\LL^*$ be the $\Z[q,q^{-1}]$-lattice spanned by the dual
canonical basis $\Bs$ of $U_q(\n)$.
It is well-known that the specialization
$\C\otimes_{\Z[q,q^{-1}]}\LL^*$,
where $\C$ is regarded as a $\Z[q,q^{-1}]$-module via $q\mapsto 1$,
is isomorphic to  the algebra $\C[N]$ of regular
functions on the group $N$.
Let $\BB^*\subset \C[N]$ denote the specialization at $q=1$
of $\B^*$ obtained in this way.

It was essentially shown by Zelevinsky \cite{Z4} that
$\RR_r$ and $\C[N]$ are isomorphic as graded algebras.
Moreover, by a dual version of Ariki's theorem \cite{A,LNT},
the basis of $\RR_r$ given by the classes of the simple modules 
is mapped to $\BB^*$ under this isomorphism.  
Therefore, for $\g$ of type $A_r$, 
the multiplication of two elements 
of respective principal degree $m$ and $n$ of $\Bs$
encodes the induction product of two simple objects
of the categories $\CC_{m,r}$ and $\CC_{n,r}$.

\subsection{}
Consider the example of \ref{A5} for type $A_5$.
The vector $b$ has principal degree $8$.
Therefore, under the above isomorphism it corresponds to 
a simple object $M$ of the category $\CC_{8,5}$.
The fact that $b$ is imaginary translates into
the fact that the induction square of $M$ is not simple.
In the Zelevinsky parametrization of simple modules by multi-segments,
the $H_8(t)$-module $M$ is labelled by the multi-segment
\[
[1,2],[2,3,4],[3],[4,5].
\]
It has dimension $252$. Its induction square $M\odot M$ is a $H_{16}(t)$-module
with two composition factors: the $2\,522\,520$-dimensional simple
module $M'$ labelled by
\[
[1,2],[2,3],[3,4],[4,5],[1,2,3,4],[2,3,4,5]
\]
and the $814\,773\,960$-dimensional simple module $M''$ labelled by
\[
[1,2],[1,2],[2,3,4],[2,3,4],[3],[3],[4,5],[4,5].
\]
 
\subsection{} Let $U_q(\slchap_N)$ be the quantized affine algebra of
type $A_{N-1}^{(1)}$ with parameter $q$ a square root of~$t$
(see for example \cite{CP} for the defining relations of $U_q(\slchap_N)$).
The quantum affine Schur-Weyl duality between $H_m(t)$ and
$U_q(\slchap_N)$ \cite{CP, Ch, GRV} gives a functor $\F_{m,N}$ from
the category $\CC_m$ to the category of level $0$ finite-dimensional
representations of $U_q(\slchap_N)$.
If $N\ge m$, $\F_{m,N}$ maps the non-zero simple modules of $\H_m(t)$
to non-zero simple modules of $U_q(\slchap_N)$. 
Moreover, 
for $M_1$ in $\CC_{m_1}$ and $M_2$ in $\CC_{m_2}$ one has
\[
\F_{m_1+m_2,N}(M_1\odot M_2) =
\F_{m_1,N}(M_1)\otimes\F_{m_2,N}(M_2)\,.
\]

Let $N\ge 16$.
The image of the $H_8(t)$-module $M$ under $\F_{8,N}$ is a
finite-dimensional irreducible representation $V$ of $U_q(\slchap_N)$
such that $V\otimes V$ is not irreducible.
The Drinfeld polynomials parametrizing $V$ in the classification
of Chari and Pressley are
\begin{eqnarray*}
P_1(u) &=& u - q^{-6},\\
P_2(u) &=& (u-q^{-3})(u-q^{-9}),\qquad\qquad\qquad\qquad \\
P_3(u) &=& u - q^{-6},\\
P_k(u) &=& 1,\quad (k\ge 4).
\end{eqnarray*}
The composition factors of $V\otimes V$ are the images $V'$ and
$V''$ of $M'$ and $M''$ under $\F_{16,N}$. 
The Drinfeld polynomials of $V'$ are
\begin{eqnarray*}
P_2(u)& =& (u-q^{-3})(u-q^{-5})(u-q^{-7})(u-q^{-9}),\\
P_4(u)& =& (u - q^{-5})(u-q^{-7}),\\
P_k(u)& =& 1,\quad (k\not= 2,4).
\end{eqnarray*}
and those of $V''$ are
\begin{eqnarray*}
P_1(u) &=& (u - q^{-6})^2,\\
P_2(u) &=& (u-q^{-3})^2(u-q^{-9})^2,\qquad\qquad\qquad\quad\\
P_3(u) &=& (u - q^{-6})^2,\\
P_k(u) &=& 1,\quad (k\ge 4).
\end{eqnarray*}
In fact, since the highest non trivial Drinfeld polynomial 
of $V''$ is $P_4(u)$, we see that $V\otimes V$ is not simple
for any $N\ge 5$.
For $N=5$, we have 
\[
\dim V = 1\,995,\qquad \dim V'=7\,350,\qquad \dim V''=3\,972\,675.
\]

\subsection{}
The meaning of Conjecture~\ref{CONJ1} in this setting should be as
follows.
\begin{conjecture}
Let $V$ and $W$ be finite-dimensional irreducible representations
of $U_q(\slchap_N)$. Suppose that $V\otimes V$ is irreducible.
Then $V\otimes W$ has simple socle $S$ and simple head $H$.
Moreover if $V\otimes W$ is not simple $S$ and $H$ are not isomorphic.
\end{conjecture}

\bigskip
\centerline{\bf Acknowledgements}

\medskip\noindent
I thank M. Kashiwara and A. Zelevinsky for fruitful discussions.
In particular, I am grateful to A. Zelevinsky who convinced me,
after I found an example of imaginary vector in type $G_2$,
that imaginary vectors should also exist for other types. 

This note was written during my visit to M.S.R.I. Berkeley
to participate in the program ``Infinite dimensional algebras
and mathematical physics'', organized by E. Frenkel,
V. Kac, I. Penkov, V. Serganova and G. Zuckerman. 
I thank the organizers for inviting me.
 

%
\bigskip
\small

\noindent
\begin{tabular}{ll}
 {\sc B. Leclerc} : &
D\'epartement de Math\'ematiques,
Universit\'e de Caen, Campus II,\\
& Bld Mar\'echal Juin,
BP 5186, 14032 Caen cedex, France\\
&email : {\tt leclerc@math.unicaen.fr}
\end{tabular}

\bigskip
\begin{thebibliography}{ABC} \scriptsize
%
\bibitem[\bf A]{A}{\sc S. Ariki}, {\it On the decomposition numbers 
of the Hecke
algebra of $G(n,1,m)$}, J. Math. Kyoto Univ. {\bf 36} (1996),
789--808.
%
\bibitem[\bf BZ]{BZ}{\sc A. Berenstein, A. Zelevinsky}, {\it String
bases for quantum groups of type $A_r$}, Adv. Soviet Mathematics,
{\bf 16} (1993), 51--89.
%
\bibitem[\bf C1]{C}{\sc P. Caldero}, {\it Adapted algebras for the
Berenstein-Zelevinsky conjecture}, Preprint 2001, 
{\tt math.RT/0104165}.
%
\bibitem[\bf C2]{C2}{\sc P. Caldero}, {\it A multiplicative 
property of quantum flag minors}, Preprint 2001,
{\tt math.RT/0112205}.
%
\bibitem[\bf CP]{CP}{\sc V. Chari, A. Pressley}, {\it Quantum affine
algebras and affine Hecke algebras}, {Pacific J. Math.} {\bf 174} (1996),
295--326.
%
\bibitem[\bf Ch]{Ch}{\sc I. V. Cherednik}, {\it A new interpretation
of Gelfand-Tzetlin bases}, {\it Duke Math. J.}, 54 (1987), 563--577.
%
\bibitem[\bf FZ]{FZ}{\sc S. Fomin, A. Zelevinsky}, {\it Cluster
algebras I: Foundations}, to appear in J. Amer. Math. Soc.,
{\tt math.RT/0104151}.
%
\bibitem[\bf GRV]{GRV}{\sc V. Ginzburg, N. Yu Reshetikhin, E. Vasserot},
{\it Quantum groups and flag varieties}, {\it A.M.S. Contemp. Math.} {\bf 175}
(1994), 101--130.
%
\bibitem[\bf G]{G}{\sc J. A. Green}, {\it Quantum groups, Hall
algebras and quantum shuffles}, in `Finite reductive groups'
(Luminy 1994), 273-290, Birkh\"auser Prog. Math. {\bf 141}, 1997.
%
\bibitem[\bf K1]{K0}{\sc M. Kashiwara}, {\it On crystal bases of
the $q$-analogue of universal enveloping algebras},
Duke Math. J. {\bf 63} (1991), 465--516.
%
\bibitem[\bf K2]{K}{\sc M. Kashiwara}, {\it Global bases of quantum groups},
Duke Math. J. {\bf 69} (1993), 455--485.          
%
\bibitem[\bf K3]{K2}{\sc M. Kashiwara}, {\it On level zero
representations of quantized affine algebras}, (2000),
 {\tt math.QA/0010293}.
%
\bibitem[\bf KS]{KS}{\sc M. Kashiwara, Y. Saito}, {\it Geometric
construction of crystal bases}, Duke Math. J. {\bf 89} (1997),
9--36.
%
\bibitem[\bf LR]{LR}{\sc P. Lalonde, A. Ram}, {\it Standard Lyndon
bases of Lie algebras and enveloping algebras},
Trans. Amer. Math. Soc. {\bf 347} (1995), 1821--1830.
%
\bibitem[\bf Le]{Le}{\sc B. Leclerc}, {\it Dual canonical bases,
quantum shuffles and $q$-characters}, Preprint 2002.
%
\bibitem[\bf LNT]{LNT}{\sc B. Leclerc, M. Nazarov, J.-Y. Thibon},
{\it Induced representations of affine Hecke algebras and 
canonical bases of quantum groups},
to appear in `Studies in memory of Issai Schur', Birkh\"auser 2002.
%
\bibitem[\bf Lo]{Lo}{\sc M. Lothaire}, {\it Combinatorics on words},
Readings, Massachusetts, 1983.
%
\bibitem[\bf L1]{L}{\sc G. Lusztig}, {\it Introduction to quantum
groups}, Birkh\"auser 1993.
%
\bibitem[\bf L2]{L2}{\sc G. Lusztig}, {\it Quivers, perverse sheaves
and quantized enveloping algebras}, J. Amer. Math. Soc. {\bf 4},
(1991), 365--421.
%
\bibitem[\bf Re]{Rei}{\sc M. Reineke}, {\it Multiplicative properties of 
dual canonical bases of quantum groups}, J. Algebra {\bf 211} (1999),
134--149.
%
\bibitem[\bf RZ]{RZ}{\sc V. Retakh, A. Zelevinsky}, {\it The base 
affine space and canonical bases in irreducible representations
of the group $Sp(4)$}, Soviet Math. Dokl. {\bf 37} (1988), 618--622.
%
\bibitem[\bf Reu]{Reu}{\sc C. Reutenauer}, {\it Free Lie algebras},
Oxford University Press, 1993.
%
\bibitem[\bf Ro1]{R1}{\sc M. Rosso}, {\it Groupes quantiques et
alg\`ebres de battage quantiques}, C. R. Acad. Sci. Paris {\bf 320}
(1995), 145-148.
%
\bibitem[\bf Ro2]{R2}{\sc M. Rosso}, {\it Quantum groups and quantum
shuffles}, Invent. Math. {\bf 133} (1998), 399--416.
%
\bibitem[\bf Ro3]{R3}{\sc M. Rosso}, {\it Lyndon bases and the multiplicative
formula for $R$-matrices}, Preprint 2002.
%
\bibitem[\bf Z1]{Z}{\sc A. Zelevinsky}, {\it Connected components of
real double Bruhat cells}, Intern. Math. Res. Notices, {\bf 21} 2000,
1131--1153.
%
\bibitem[\bf Z2]{Z2}{\sc A. Zelevinsky}, {\it Personal communication}.
%
\bibitem[\bf Z3]{Z3}{\sc A. Zelevinsky}, {\it From Littlewood-Richardson
coefficients to cluster algebras in three lectures},
{\tt math.RT/0112062}.
%
\bibitem[\bf Z4]{Z4}{\sc A. Zelevinsky}, {\it Induced
representations of reductive $p$-adic groups II}, 
Ann. Sci. E.N.S. {\bf 13} (1980), 165--210.
%
\end{thebibliography}
\end{document}